\documentclass[12pt]{amsart} 
\usepackage{amsfonts,graphics,amsmath,amsthm,amsfonts,amscd, amssymb,amsmath,latexsym,multicol}
\usepackage{epsfig}
\usepackage{flafter}

\input diagrams

\makeatletter

\def\jobis#1{FF\fi
  \def\predicate{#1}%
  \edef\predicate{\expandafter\strip@prefix\meaning\predicate}%
  \edef\job{\jobname}%
  \ifx\job\predicate
}

\makeatother

\if\jobis{proposal}%
 \def\try{subsection}%
\else
  \def\try{section}%
\fi

\theoremstyle{plain}
\newtheorem{theorem}{Theorem}[\try]
\newtheorem{corollary}[theorem]{Corollary}
\newtheorem{lemma}[theorem]{Lemma}
\newtheorem{claim}[theorem]{Claim}

\newtheorem{proposition}[theorem]{Proposition}
\newtheorem{definition-lemma}[theorem]{Definition-Lemma}
\newtheorem{question}[theorem]{Question}
\newtheorem{definition}[theorem]{Definition}
\newtheorem{remark}[theorem]{Remark}

\newtheorem{conjecture}[theorem]{Conjecture}

\newtheorem{problem}[theorem]{Problem}

\def\ideal#1.{I_{#1}}
\def\ring#1.{\mathcal {O}_{#1}}
\def\fring#1.{\hat{\mathcal {O}}_{#1}}
\def\proj#1.{\mathbb {P}(#1)}
\def\pr #1.{\mathbb {P}^{#1}}
\def\dpr #1.{\hat{\mathbb {P}}^{#1}}
\def\af #1.{\mathbb A^{#1}}
\def\Hz #1.{\mathbb F_{#1}}
\def\Hbz #1.{\overline{\mathbb F}_{#1}}
\def\fb#1.{\underset #1 {\times}}
\def\rest#1.{\underset {\ \ring #1.} \to \otimes}
\def\au#1.{\operatorname {Aut}\,(#1)}
\def\deg#1.{\operatorname {deg } (#1)}
\def\pic#1.{\operatorname {Pic}\,(#1)}
\def\pico#1.{\operatorname{Pic}^0(#1)}
\def\picg#1.{\operatorname {Pic}^G(#1)}
\def\ner#1.{NS (#1)}
\def\rdown#1.{\llcorner#1\lrcorner}
\def\rfdown#1.{\lfloor{#1}\rfloor}
\def\rup#1.{\ulcorner{#1}\urcorner}
\def\rcup#1.{\lceil{#1}\rceil}
\def\cone#1.{\operatorname {NE}(#1)}
\def\ccone#1.{\overline{\operatorname {NE}}(#1)}
\def\coef#1.{\frac{(#1-1)}{#1}}
\def\vit#1.{D_{\langle #1 \rangle}}
\def\mm#1.{\overline {M}_{0,#1}}
\def\H1#1.{H^1(#1,{\ring #1.})}
\def\ac#1.{\overline {\mathbb F}_{#1}}

\def\adj#1.{\frac {#1-1}{#1}}
\def\spn#1.{\overline{#1}}
\def\pek#1.#2.{\Cal P^{#1}(#2)}
\def\plk#1.#2.{\Cal P^{\leq #1}(#2)}
\def\ev#1.{\operatorname{ev_{#1}}}
\def\ilist#1.{{#1}_1,{#1}_2,\dots}
\def\bminv#1.{(\nu_1,s_1;\nu_2,s_2;\dots ;\nu_{#1},s_{#1};\nu_{r+1})}
\def\zinv#1.{(\nu_1,s_1;\nu_2,s_2;\dots ;\nu_{#1},s_{#1};0)}
\def\iinv#1.{(\nu_1,s_1;\nu_2,s_2;\dots ;\nu_{#1},s_{#1};\infty)}

\def\llist#1.#2.{{#1}_1,{#1}_2,\dots,{#1}_{#2}}
\def\lomitlist#1.#2.{{#1}_1,{#1}_2,\dots,\hat {{#1}_i}, \dots, {#1}_{#2}}
\def\lomitlistz#1.#2.{{#1}_0,{#1}_1,\dots,\hat {{#1}_i}, \dots, {#1}_{#2}}
\def\loc#1.#2.{\Cal O_{#1,#2}}
\def\fderiv#1.#2.{\frac {\partial #1}{\partial #2}}
\def\deriv#1.#2.{\frac {d #1}{d #2}}
\def\map#1.#2.{#1 \longrightarrow #2}
\def\rmap#1.#2.{#1 \dasharrow #2}
\def\emb#1.#2.{#1 \hookrightarrow #2}
\def\non#1.#2.{\text {Spec }#1[\epsilon]/(\epsilon)^{#2}}
\def\Hi#1.#2.{\text {Hilb}^{#1}(#2)}
\def\sym#1.#2.{\operatorname {Sym}^{#1}(#2)}
\def\Hb#1.#2.{\text {Hilb}_{#1}(#2)}
\def\Hm#1.#2.{\Hom_{#1}(#2)}
\def\prd#1.#2.{{#1}_1\cdot {#1}_2\cdots {#1}_{#2}}
\def\Bl #1.#2.{\operatorname {Bl}_{#1}#2}
\def\pl #1.#2.{#1^{\otimes #2}}
\def\mgn#1.#2.{\overline {M}_{#1,#2}}
\def\ialist#1.#2.{{#1}_1 #2 {#1}_2, #2\dots}
\def\pair#1.#2.{\langle #1, #2\rangle}
\def\vandermonde#1.#2.{\left|
\begin{matrix}
1 & 1 & 1 & \dots & 1\\
{#1}_1 & {#1}_2 & {#1}_3 & \dots & {#1}_{#2}\\
{#1}_1^2 & {#1}_2^2 & {#1}_3^2 & \dots & {#1}_{#2}^2\\
\vdots & \vdots & \vdots & \ddots & \vdots\\
{#1}_1^{#2-1} & {#1}_2^{#2-1} & {#1}_2^{#2-1} & \dots & {#1}_{#2}^{#2-1}\\
\end{matrix}
\right|
}
\def\vandermondet#1.#2.{\left|
\begin{matrix}
1 & {#1}_1   & {#1}_1^2 & \dots & {#1}_1^{#2-1}\\
1 & {#1}_2   & {#1}_2^2 & \dots & {#1}_2^{#2-1}\\
1 & {#1}_3   & {#1}_3^2 & \dots & {#1}_3^{#2-1}\\
\vdots & \vdots & \vdots & \ddots & \vdots\\
1 & {#1}_{#2}& {#1}_{#2}^2 & \dots & {#1}_{#2}^{#2-1}\\
\end{matrix}
\right|
}
\def\gr#1.#2.{\mathbb{G}(#1,#2)}

\def\alist#1.#2.#3.{{#1}_1 #2 {#1}_2 #2\dots #2 {#1}_{#3}}
\def\zlist#1.#2.#3.{#1_0 #2 #1_1 #2\dots #2 #1_{#3}}
\def\lomitlist30#1.#2.#3.{{#1}_0,{#1}_1 #2 \dots #2\hat {{#1}_i} #2\dots #2 {#1}_{#3}}
\def\lmap#1.#2.#3.{#1 \overset{#2}{\longrightarrow} #3}
\def\mes#1.#2.#3.{#1 \longrightarrow #2 \longrightarrow #3}
\def\ses#1.#2.#3.{0\longrightarrow #1 \longrightarrow #2 \longrightarrow #3 \longrightarrow 0}
\def\les#1.#2.#3.{0\longrightarrow #1 \longrightarrow #2 \longrightarrow #3}
\def\res#1.#2.#3.{#1 \longrightarrow #2 \longrightarrow #3\longrightarrow 0}
\def\Hi#1.#2.#3.{\text {Hilb}^{#1}_{#2}(#3)}
\def\ten#1.#2.#3.{#1\underset {#2}{\otimes} #3}
\def\lomitlist30#1.#2.#3.{{#1}_0 #2 {#1}_1 #2 \dots #2 \hat {{#1}_i} #2 \dots #2 {#1}_{#3}}

\def\Hom{\operatorname{Hom}}

\def\sp{\operatorname{Spec}}

\def\dim{\operatorname{dim}}

\def\deg{\operatorname{deg}}

\def\im{\operatorname{Im}}

\def\lcc{\operatorname{LLC}}
\def\lcs{\operatorname{LCS}}

\def\Div{\operatorname{Div}}

\def\mult{\operatorname{mult}}

\def\rest{\operatorname{res}}

\def\vol{\operatorname{vol}}

\def\p{\mathbb P}

\def\e{\Cal E}

\def\e1{E_1}
\def\e2{E_2}

\def\qle{\sim_{\mathbb Q}}

\def\mapdown#1{\big\downarrow\rlap{$\vcenter
{\hbox{$\scriptstyle#1$}}$}}

\def\mapse#1{
{\vcenter{\hbox{$\mathop{\smash{\raise1pt\hbox{$\diagdown$}\!\lower7pt
\hbox{$\searrow$}}\vphantom{p}}\limits_{#1}\vphantom{\mapdown{}}$}}}}

\def\VR#1.{height#1pt&\omit&&\omit&&\omit&&\omit&&\omit&\cr}

\def\VRT#1.{height#1pt&\omit&&\omit&\cr}

\begin{document}
\title{Boundedness of pluricanonical maps of varieties of general type}
\author{Christopher D. Hacon} 
\address{Department of Mathematics \\  
University of Utah\\  
155 South 1400 E\\
JWB 233\\
Salt Lake City, UT 84112, USA}
\email{hacon@math.utah.edu}
\author{James M\textsuperscript{c}Kernan} 
\address{Department of Mathematics\\ 
University of California at Santa Barbara\\ 
Santa Barbara, CA 93106, USA} 
\email{mckernan@math.ucsb.edu}

\thanks{The first author was partially supported by NSA research grant no:
  MDA904-03-1-0101 and by a grant from the Sloan Foundation.  We are in debt to
  \cite{Tsuji99} and the beautiful ideas and techniques introduced in that paper.  We
  would like to thank Meng Chen for pointing out that a previous version of
  \eqref{q_three} was too optimistic.  We would also like to thank the referee, whose
  extensive comments were very useful.}

\begin{abstract} Using the techniques of \cite{Siu98} and \cite{Kawamata99}, we prove that
certain log forms may be lifted from a divisor to the ambient variety.  As a consequence
of this result, following \cite{Tsuji99}, we show that: For any positive integer $n$ there
exists an integer $r_n$ such that if $X$ is a smooth projective variety of general type
and dimension $n$, then $\phi _{rK_X}\colon \rmap X.{\proj {H^0(\ring X.(rK_X))}.}.$ is
birational for all $r\geq r_n$.
\end{abstract}
\maketitle
\section{Introduction}

One of the main problems of complex projective algebraic geometry is to understand the
structure of the pluricanonical maps.  If $X$ is of general type, then by definition the
pluricanonical map $\phi _{rK_X}\colon \rmap X.{\proj {H^0(\ring X.(rK_X))}.}.$ is
birational for all sufficiently large $r$.  It is a natural question to then ask if there
is an $r_n$ such that $\phi_{rK_X}$ is birational as soon as $r\geq r_n$, uniformly for
any variety of general type of dimension $n$.  When $X$ is a curve, it is well known that
$\phi_{rK_X}$ is birational for $r\geq 3$ and when $X$ is a surface, Bombieri proved in
\cite{Bombieri70} that $\phi_{rK_X}$ is birational, for $r\geq 5$.

However, starting with threefolds the problem is substantially harder, since there are
threefolds of general type for which the minimal model is necessarily singular.  In fact
it is easy to write down examples of threefolds of general type, for which the index of
$K_X$ (the smallest multiple of $K_X$ which is Cartier) on a minimal model is arbitrarily
large (for example, take any hypersurface of sufficiently large degree in any weighted
projective space).  In this case, a priori, the degree of $K_X$ (equivalently the volume
of $K_X$ on a smooth model) could be arbitrarily small.  Even the case when the minimal
model is smooth, or at least the index of $K_X$ is bounded, has attracted considerable
attention, see for example, \cite{Hanamura85}, \cite{Benveniste86}, \cite{Matsuki86} and
\cite{Cheng01}, and there have only been partial results in the general case, see
\cite{Kollar86}, \cite{Luo94} and \cite{Luo00}.

Recently, Tsuji settled this problem for $3$-folds and assuming the minimal model program,
for all higher dimensional varieties, see \cite{Tsuji99}.  Using ideas of Tsuji, we prove
the following:

\begin{theorem}\label{t_birational} For any positive integer $n$, there exists an integer $r_n$ 
such that if $X$ is a smooth projective variety of general type and dimension $n$, then
$\phi_{rK_X}\colon \rmap X.{\proj {H^0(\ring X.(rK_X))}.}.$ is birational for all $r\geq r_n$.
\end{theorem}

\eqref{t_birational} has some interesting consequences:
\begin{corollary}\label{c_1} 
Let $n$ be a positive integer and $M$ be a positive constant. Then the family of all
smooth projective varieties of general type of dimension $n$ and volume of $K_X$ at most
$M$, is birationally bounded.
\end{corollary}

\begin{corollary}\label{c_2} Let $n$ be a positive integer. Then there is a 
positive constant $\eta _n>0$ such that given any smooth projective
variety $X$ of general type of dimension $n$, the volume of $K_X$
is at least $\eta _n$.
\end{corollary}
As Maehara points out in \cite{Maehara83} (see also \cite{Tsuji99}), one also recovers the
Severi-Iitaka conjecture:
\begin{corollary}\label{c_3} 
For any fixed variety $X$, there exist only finitely many dominant rational maps $\pi
\colon \rmap X.Y.$ (modulo birational equivalence), where $Y$ is a variety of general
type.
\end{corollary}

The proof of the above results closely follows Tsuji's original strategy; in fact this
note grew out of our desire to understand \cite{Tsuji99}.  Even though our proof is based
on the ideas of this paper, it only relies on the techniques of algebraic geometry.  In
particular we do not use the ``Analytic Zariski Decomposition'' and we do not make use of
the results of \cite{Tsuji01} (which we were unable to follow).  The main tools employed
in this paper are the algebraic techniques of the minimal model program and in particular
the theory of log canonical centres and of multiplier ideals.

We now give a short informal sketch of the proof of \eqref{t_birational}, which closely
follows the strategy of Tsuji's proof, see \cite{Tsuji99}.  As one would expect, the goal
is to produce a $\mathbb{Q}$-divisor $\Delta\sim _{ \mathbb{Q}}\lambda K_X$ such that
$\Delta$ has an isolated log canonical centre at a general point $x\in X$. There is a well
established strategy for producing divisors with log canonical centres at a general point
$x\in X$, as soon as the volume is sufficiently large, so that the rational number
$\lambda$ depends only on the volume of $K_X$.

The problem is that the smallest log canonical centre $V$ containing $x$ might well be of
positive dimension.  In order to produce an isolated log canonical centre, one proceeds by
inductively cutting down the (dimension of the) log canonical centres.  As $x$ is general,
a resolution $W$ of $V$ must have general type, and so by induction we can certainly find 
a $\mathbb{Q}$-divisor $\Theta$ on $W$ with the required properties.  

At this point, the idea is to lift sections from $V$ to the whole of $X$.  In principle
this ought to be straightforward; indeed if $K_X$ is nef, then we can lift $\Theta$ to the
whole of $X$, as an easy application of Serre vanishing.  In practice we cannot assume
that $K_X$ is nef, and this step of the argument forms the technical heart of the paper.
In fact we generalise some results and techniques of \cite{Siu98} and \cite{Kawamata99}
for extending pluricanonical forms from a divisor to the whole space.  Instead of just
lifting pluricanonical forms, we lift certain log pluricanonical forms, from a log
canonical centre to the whole of $X$, using a variant of multiplier ideal sheaves, see \S
\ref{s_lift}.  It is then not too hard to prove that we can then lift log canonical
centres from $W$ to $X$, see \S \ref{s_centres}.

Then we prove, assuming \eqref{t_birational} in dimension less than $n$, that
$\phi_{rK_X}$ is birational, where $r$ is a fixed linear function of $1/\vol(X)^{1/n}$,
which only depends on the dimension.  If $\vol (K_X)\geq 1$, the theorem follows easily.
If $\vol (K_X)<1$, one then shows that $X$ belongs to a birationally bounded family and
the existence of $r_n$ is then clear.

As already pointed out, the main technical point is \eqref{c_ext}, the result which shows
that certain log forms may be lifted from a divisor to the ambient variety.  We hope that
it will find other applications in a variety of contexts. For example, in \cite{HM05a}, we
use this result to prove a conjecture of Shokurov which in particular implies that the
fibers of a resolution of a variety with log terminal singularities are rationally chain
connected.

We now raise some questions and problems, which are closely related to
\eqref{t_birational}:

\begin{problem}\label{p_explicit} Find explicit (hopefully small) values of $r_n$.  
\end{problem}

Even determining $r_3$ would seem interesting.  Indeed in the past twenty years or so, our
knowledge of threefolds has increased considerably; in particular we have Reid's powerful
Riemann-Roch formula for terminal threefolds and we are able to write down many explicit
examples of threefolds of general type.  Thus determining the value of $r_3$ would seem to
be a reasonable problem.  In fact we raise the following:

\begin{question}\label{q_three} Is the rational map $\phi_{27K_X}$ always birational,
whenever $X$ is a threefold of general type?
\end{question}

Note that we have little or no evidence for an affirmative answer to \eqref{q_three},
apart from the examples of Fletcher in \cite{Fletcher00}.  In particular for the threefold
$X_{46}$, a hypersurface in the weighted projective space $\proj 4,5,6,7,23.$, see
\cite{Fletcher00}, $\phi_{mK_X}$ is birational iff $m=23$, or $m\geq 27$.

Finally it seems natural to conjecture that we can drop the condition that $X$ is of
general type in \eqref{t_birational}:

\begin{conjecture}\label{c_general} There is a positive integer $r_n$ with the following 
property:  

Let $X$ be a smooth variety of non-negative Kodaira dimension, of dimension $n$.  Then the
rational map $\phi_{rK_X}$ is birationally equivalent to the Iitaka fibration.
\end{conjecture}

 Unfortunately the methods of this paper don't seem to apply directly to \eqref{c_general}. 
\section{Preliminaries}
\label{s_pre}

\subsection{Notation and conventions}

We work over the field of complex numbers $\mathbb{C}$.  A $\mathbb{Q}$-Cartier divisor
$D$ on a normal variety $X$ is nef if $D\cdot C\geq 0$ for any curve $C\subset X$. We say
that two $\mathbb{Q}$-divisors $D_1,D_2$ are $\mathbb{Q}$-linearly equivalent ($D_1\sim
_{\mathbb{Q}} D_2$) if there exists an integer $m>0$ such that $mD_i$ are linearly
equivalent.  Given a morphism of normal varieties $f\colon \map X.Y.$, we say that two
divisors $D_1$ and $D_2$ on $X$ are $f$-linearly equivalent ($D_1\sim _fD_2$) if there is
a Cartier divisor $B$ on $Y$ such that $D_1\sim D_2+f^*B$.  We say that $D_1$ and $D_2$
are $f$-numerically equivalent ($D_1\equiv _fD_2$) if there is a $\mathbb{Q}$-Cartier
divisor $B$ on $Y$ such that $D_1\equiv D_2+f^*B$.  A {\bf pair} $(X,\Delta )$ is a normal
variety $X$ and a $\mathbb{Q}$-Weil divisor $\Delta$ such that $K_X+\Delta$ is
$\mathbb{Q}$-Cartier.  We will say that a pair $(X,\Delta)$ is a log pair, if in addition
$\Delta$ is effective.  We say that $(X,\Delta)$ is a {\bf smooth pair} if $X$ is smooth
and $\Delta$ is a $\mathbb{Q}$-divisor with simple normal crossings support.  A projective
morphism $\mu \colon \map Y.X.$ is a {\bf log resolution} of the pair $(X,\Delta )$ if $Y$
is smooth and $\mu ^{-1}(\Delta )\cup \{\text{exceptional set of } \mu \}$ is a divisor
with simple normal crossings support.  We write $\mu ^* (K_X +\Delta )=K_Y +\Gamma$ and
$\Gamma =\sum a _i\Gamma _i$ where $\Gamma _i$ are distinct reduced irreducible divisors.
The pair $(X,\Delta )$ is {\bf kawamata log terminal} (resp. {\bf log canonical}) if there
is a log resolution $\mu \colon \map Y.X.$ as above such that the coefficients of $\Gamma$
are strictly less than one i.e. $a _i<1$ for all $i$ (resp.  $a _i\leq 1$).  The number
$1-a _i$ is the {\bf log discrepancy} of $\Gamma _i$ with respect to $(X, \Delta )$.  We
say that a subvariety $V\subset X$ is a {\bf log canonical centre} if it is the image of a
divisor of log discrepancy at most zero.  We will denote by $\lcc (X,\Delta ,x)$ the set
of all log canonical centres containing a point $x\in X$.  A {\bf log canonical place} is
a valuation corresponding to a divisor of log discrepancy at most zero.  A log canonical
centre is {\bf pure} if $K_X+\Delta$ is log canonical at the generic point of $V$.  If
moreover there is a unique log canonical place lying over the generic point of $V$, we
will say that $V$ is an {\bf exceptional log canonical centre}.  

\subsection{Volumes}
\begin{definition}\label{d_volume} Let $X$ be an irreducible
projective variety and 
let $D$ be a big
$\mathbb{Q}$-divisor.  The \textbf{volume of $D$}, 
is 
$$
\vol(D)=\limsup _{m\to\infty}\frac {n!h^0(X,mD)}{m^n}.
$$
\end{definition}
When $D$ is very ample, $\vol (D)$ is just the degree of the image of $X$ 
in $\p H^0(X,\ring X.(D))$. When $D$ is nef, then $\vol (D)=D^{\dim (X)}$.
It turns out that the volume of $D$ only depends on its numerical class
and one can extend the volume function to 
a continuous function $\vol \map \colon N^1_\mathbb{R}(X).\mathbb{R}.$ (cf. \cite{Lazarsfeld04b}).

\begin{lemma}\label{l_low}
Let $X$ be a projective variety, 
$D$ a divisor such that $\phi _{D}$ is birational with image $Z$. Then, 
the volume of $D$ is at least the degree of $Z$ and hence at least $1$.
\end{lemma}
\begin{proof} We may replace $X$ by an appropriate birational model
and in particular we may assume that $\phi :=\phi _{D}$ is a morphism.
Since $\phi $ is birational, its image $Z$ is a non-degenerate subvariety
of projective space $\pr N.$ and hence its degree is at least $1$.
From the inclusion 
$$
\emb {\phi ^* \ring {\pr N.}.(1)|Z}.{\ring X.(D)}.,
$$
one sees that $\vol (D)\geq \vol (\ring {\pr N.}.(1)|Z) \geq 1$.
\end{proof}

\subsection{Log Canonical Centres}

Here we recall some well known results regarding log canonical centres.  The first Lemma
shows that when the volume of a $\mathbb{Q}$-divisor $L$ is big, it is possible to produce a
numerically equivalent $\mathbb{Q}$-divisor $\Delta _x$ with high multiplicity at a general point
$x$.

\begin{lemma}\label{l_sing} Let $V$ be an irreducible projective variety of dimension $d$, $L$ a 
big $\mathbb{Q}$-Cartier divisor on $V$, and $x\in V$ a smooth point. If for a positive rational
number $\alpha$, one has $$\vol (L)>\alpha ^d,$$ then for any sufficiently divisible
integer $k\gg 0$, there exists a divisor
$$
A=A_x \in |kL| \qquad \text{with} \qquad \mult _x(A)>k\alpha .
$$
\end{lemma}
This is easily seen by a parameter count comparing the number of sections in $\ring V.(kL)$
and the number of conditions one needs to impose for a function to vanish at a smooth
point to order at least $k\alpha$.  The next lemma shows that if the multiplicity of a
$\mathbb{Q}$-divisor at a given point is sufficiently big, then the given point belongs to an
appropriate log canonical centre.

\begin{lemma}\label{l_d} Let $(X,\Delta )$ be a log pair, $x$ a smooth point of
$X$. If $\mult _x(\Delta )\geq \dim X$, then $\lcc (X,\Delta ,x)\ne \emptyset$.  If $\mult
_x(\Delta )<1$, then $\lcc (X,\Delta ,x){=} \emptyset$.
\end{lemma}
\begin{proof} See \cite{Lazarsfeld04b} 9.3.2 and 9.5.13.
\end{proof}
Assume that $(X,\Delta )$ is
log canonical at $x$.
It is often useful to assume that $Z=\lcc (X,\Delta ,x)$
is an irreducible variety. To this end recall (cf. \cite{Kawamata97}
and \cite{Ambro98} (3.4))
\begin{lemma}\label{l_irr}
Let $X$ be a normal variety and let $\Delta$ be an effective $\mathbb{Q}$-Cartier divisor
such that $K_X+\Delta$ is $\mathbb{Q}$-Cartier. Assume that $x\in X$ is a kawamata log
terminal point of $X$ and that $(X,\Delta )$ is log canonical near $x\in X$. If
$W_1,W_2\in \lcc (X,\Delta , x)$, and $W$ is an irreducible component of $W_1\cap W_2$
containing $x$, then $W$ is in $\lcc (X, \Delta , x)$.  Therefore, if $(X,\Delta )$ is not
kawamata log terminal at $x$, then $\lcc(X,\Delta , x)$ has a unique minimal irreducible
element, say $V$.  Moreover, there exists an effective $\mathbb{Q}$-divisor $E$ such that
$$
\lcc (X,(1-\epsilon) \Delta +\epsilon E ,x)=\{ V\}
$$
for all $0< \epsilon \ll 1$. We may also assume that there is a unique log canonical place
lying above $V$.  If $X$ is projective, $x\in X$ is general and $L$ is a big divisor, then
one can take $E=aL$ for some positive number $a$.
\end{lemma}

\subsection{Non-vanishing and birational maps}

Recall the following well known method for producing sections of adjoint line bundles:

\begin{lemma}\label{l_lift} Let $X$ be a  smooth projective variety and
$D$ a big divisor on $X$.  Let $0<\lambda <1$ a rational number and assume that for every
general point $x\in X$, there exists $\Delta=\Delta_x\sim _{\mathbb{Q}} \lambda D$ a
$\mathbb{Q}$-divisor such that $x$ is a maximal (in terms of inclusion) element of $\lcs
(X,\Delta,x)$.  Then $h^0(\mathcal {O}_X(K_X+D))>0$ and $h^0(\mathcal {O}_X(K_X+2D))\geq
2$.
\end{lemma}
\begin{proof} For $m\gg 0$, let $G$ be a general divisor in $|mD|$.  Then $x$ is not
contained in the support of $G$ and so $x$ is the only element of $\lcc(X,\Delta
+\frac{1-\lambda}{m}G, x)$.  Let $\mu \colon \map X'.X.$ be a log resolution of $(X,\Delta
+\frac{1-\lambda}{m}G)$ and let
$$
\mathcal{I}:=\mathcal{I}(X,\Delta+\frac{1-\lambda}{m}G):= \mu_*\ring X'.(K_{X'/X}-\rdown \mu ^* (\Delta
+\frac{1-\lambda}{m}G).) .
$$ 
We may assume that $\mu ^* |mD|=F+|M|$ where $F$ has simple normal crossings and $|M|$ is
base point free and hence $M$ is nef and big.  Then
$$
\mu^*D-\rdown\mu ^* (\Delta +\frac{1-\lambda}{m}G).\sim _{\mathbb{Q}} \mu^*D-\rdown\mu ^*
(\Delta +\frac{1-\lambda}{m}F).\sim _{\mathbb{Q}}
$$
$$
\{\mu ^* (\Delta +\frac{1-\lambda}{m}F)\}+\frac{1-\lambda}{m}M .
$$
So $\mu ^*D-\rdown\mu ^* (\Delta +\frac{1-\lambda}{m}G).$ is numerically equivalent to the
sum of a fractional divisor with simple normal crossings and a nef and big divisor. By
Kawamata-Viehweg vanishing, one sees that $R^i\mu _*\ring X'.(K_{X'}+\mu^*D-\rdown \mu ^*
(\Delta +\frac{1-\lambda}{m}G).)=0$ for $i>0$.  So there is a short exact sequence
$$
\ses {\ring X.(K_X+D)\otimes \mathcal{I}}.{\ring X.(K_X+D)}.{\frac{\ring X.(K_X+D)}{\ring X.(K_X+D)\otimes \mathcal{I}}}..
$$ 
Again by Kawamata-Viehweg vanishing, one has that
$$
h^1(\ring X.(K_X+D)\otimes \mathcal{I})=
h^1(\ring X'.(K_{X'}+\mu^*D-\rdown\mu ^* (\Delta +\frac{1-\lambda}{m}G).))=0.
$$
Since $x$ is the only element of $\lcc (X,\Delta +\frac{1-\lambda}{m}G, x)$, one has
that $\mathbb{C} _x$ is a direct summand of $\frac{\ring X.(K_X+D)}{\ring X.  (K_X+D)\otimes
  \mathcal{I}}$ and so $\ring X. (K_X+D)$ is generated at $x$.  

Pick a general point $x_1$.  Then we may find $D_1\qle \lambda D$ with an isolated log
canonical centre at $x_1$.  Now pick a general point $x_2$, not in the support of $D_1$,
and pick a general divisor $D_2\qle \lambda D$ with an isolated log canonical centre at
$x_2$.  As $x_2$ and $D_2$ are general, $D_2$ does not contain $x_1$.  Thus $x_1$ and
$x_2$ are isolated log canonical centres of $D_1+D_2\qle 2\lambda D$.

 The rest of the proof that $h^0(\ring X. (K_X+2D))\geq 2$ now follows the line above.
\end{proof}

Given a line bundle with at least two sections, one has a rational map to $\p ^1$.  The
following Lemma gives a way to produce birational maps.

\begin{lemma}\label{l_bir} Let $X$ be a smooth projective variety. 
Suppose that $h^0(X,\ring X. (mK_X))\geq 2$ for all $m\geq m_0$ and let $\map X'.{\pr
  1.}.$ be any morphism induced by sections of $\ring X. (m_0K_X)$ on an appropriate
birational model $X'$ of $X$. If $F$ denotes the general fiber and $|s K_F|$ induces a
birational morphism, then $|tK_X|$ induces a birational morphism for all $t\geq
m_0(2s+2)+s$.
\end{lemma}
\begin{proof} Following Theorem 4.6 of \cite{Kollar86} and its proof,
one sees that $|(m_0(2s+1)+s)K_X|$ gives a birational map.
Since $mK_X$ is effective for all $m\geq m_0$, the assertion follows.
\end{proof}

 It is also possible to prove a converse to \eqref{l_lift}:
\begin{lemma}\label{l_birational} Let $(X,\Delta)$ be a log pair of dimension $n$, and let 
$D$ be an integral Weil divisor such that the image $Y$ of the rational map $\phi_D$ 
has dimension $n$.  

Then there is an open set $U$ such that for any point $x\in U$, we may find $\Delta'\sim
nD$, such that $x$ is the only element of $\lcc(X,\Delta+\Delta',x)$.
\end{lemma}
\begin{proof} Let $U$ be the complement of the singular locus, the support of 
$\Delta$, and the locus where the map $\phi_D$ is not an \'etale morphism.  

 Pick $x\in X$ and let $y\in Y$ be a point of the image $\phi_D\colon\map X.Y.$.  
By assumption, $Y$ sits inside projective space $\pr N.$.  Pick $n+1$ general hyperplanes
which contain $y$, $\llist H.n+1.$ and $n+1$ rational numbers $\llist a.n+1.$, such 
that 
$$
0<a_i<1 \qquad \text{and } \qquad \sum _i a_i=n.
$$
Set 
$$
\Gamma=\sum a_iH_i.
$$
Then $y$ is the only log canonical centre of $K_Y+\Gamma$.  Let $W$ be the normalisation 
of the graph of $\phi=\phi_D$, so that there is a commutative diagram
$$
\begin{diagram}
 &   &   W    &   & \\
&\ldTo^p &  & \rdTo^q &\\
 X  & &   \rDashto^{\phi}    & & Y.
\end{diagram}
$$

Let $D_i=p_*q^*H_i$, and set 
$$
\Delta'=\sum a_iD_i.
$$
Then each $D_i$ is linearly equivalent to $D$, so that $\Delta'$ is linearly equivalent
to $nD$. On the other hand, as $\phi$ is \'etale over $y$, it follows that $x$ is an
isolated log canonical centre of $K_X+\Delta+\Delta'$.
\end{proof}

\subsection{Log additivity of the Kodaira dimension}
\label{s_additivity}

\begin{lemma}\label{l_wp} Let $X$ and $Y$ be smooth projective varieties and 
$\Delta$ an effective $\mathbb{Q}$-divisor on $X$ with simple normal crossings support.
Let $\pi\colon\map X.Y.$ be a morphism and suppose that $K_X+\Delta$ is log
canonical on the general fiber $W$ of $\pi$, and $\kappa(W,(K_X+\Delta)|_W)\geq 0$. Then
$$
\pi_*\ring X.(m(K_{X/Y}+\Delta )),
$$
is weakly positive, for any $m$ which is sufficiently divisible. 
In particular for any ample line bundle $H$ on $Y$ and
any rational number $\epsilon>0$, one has that
$\kappa (K_{X/Y}+\Delta +\epsilon \pi ^* H)\geq \dim Y$.
\end{lemma}
\begin{proof} It suffices to prove that
$$
\pi_*\ring X.(m(K_{X/Y}+\Delta_h)),
$$
is weakly positive, where $\Delta_h$ consists of those components of $\Delta$ which
dominate $Y$.  In fact both \cite{Campana01}, (4.13) and \cite{Lu02}, (9.8) prove that
this sheaf is weakly positive, as soon as $m\Delta_h$ is integral. \end{proof}

\begin{remark}\label{r_additivity} In fact both \cite{Campana01} and \cite{Lu02} prove 
a much stronger result, by putting an orbifold structure on the base which takes into
account the multiple fibers of the morphism $\pi$.  We will not need this stronger form.  
\end{remark}

The following consequences of \eqref{l_wp} are well known:

\begin{corollary}\label{c_nm} Let $X,Y, \Delta ,H $ be as in \eqref{l_wp}. 
\begin{enumerate}
\item If $K_Y$ is pseudoeffective, then for all rational numbers $\epsilon >0$
$$\kappa (K_X+\Delta +\epsilon \pi ^* H)\geq \dim Y.$$
\item If $Y$ is of general type, $\kappa (K_X+\Delta )\geq \dim Y$.
\end{enumerate}
\end{corollary}

\begin{corollary}\label{c_wp} Let $X$, $Y$, and $\Delta $ be as in \eqref{l_wp}.

Then there exists an ample line bundle $A$ on $X$ such that for all $m$ sufficiently big
and divisible, $h^0(X,\ring X.(m(K_{X/Y}+\Delta))+A)\geq 0$.
\end{corollary}
\begin{proof} Let $n=\dim X$.  Fix a very ample line bundle $H$ on $X$ such that
$A=K_X+(n+2)H$ is also ample.  By \eqref{l_wp}, for all $m$ sufficiently big and
divisible, $\kappa (m(K_{X/Y}+\Delta)+H)\geq 0$.  By \cite{Lazarsfeld04b} 11.2.13, one sees
that $h^0(X,\ring X.(K_X+(n+2)H+m(K_{X/Y}+\Delta)))\geq 0$.
\end{proof}


\section{Lifting sections}
\label{s_lift}

In this section, we prove that certain log pluricanonical forms extend from a divisor to
the ambient space.  We follow the exposition of \cite{Kawamata99} very closely, which
proves that pluricanonical forms extend, indicating those places where the definitions and
proofs have to be modified.  

\begin{definition}\label{d_multiplier} Let $(X,\Delta)$ be a smooth log pair, where every 
component of $\Delta$ has coefficient one and let $D$ be a divisor on $X$.  Let
$\mu\colon\map W.X.$ be a birational morphism from a smooth variety $W$ and let $\Theta$
be the support of all divisors of log discrepancy zero.  Assume that there is a
decomposition
$$
\mu^*D=P+M,
$$
in $\Div(W)\otimes\mathbb{R}$ such that $P$ is $\mu$-nef and $M$ is effective, where
$\Theta\cup M\cup \{$the exceptional set of  $\mu\}$ has normal crossing support, and $\Theta$ and $M$ have no components in common.
The \textbf{multiplier ideal sheaf} $\mathcal{I}_{\Delta,M}$ is defined by the following
formula
$$\mathcal{I}_{\Delta,M}=\mu_*\ring W. (K_{W/X}+\Theta -\mu ^*\Delta - \rdown M . ).$$
\end{definition}
Note that any component of the strict transform of $\Delta$ is automatically a component
of $\Theta$.  
Equivalently, in the notation analogous to \cite{Kawamata99}, one has that 
$$ \mu_*\ring W.(\rup P.+K_W+\Theta)=\mathcal{I}_{\Delta,M}(D+K_X+\Delta).$$
When $\Delta =0$, this reduces to the usual definition of $\mathcal{I}_{M}$.
The notation $\rup P.$ is somewhat misleading as $\rup P. \ne \rup P'.$
for two $\mathbb{Q}$-linearly equivalent divisors. Therefore we replace it by
$\mu ^* D-\rdown M . $.

\begin{remark}\label{r_multiplier} The definition of $\Theta$ is motivated by two 
conflicting requirements:
\begin{enumerate} 
\item On the one hand, we require 
$$
(K_W+\Theta)-\mu^*(K_X+\Delta),
$$
to be effective and exceptional.  This forces us to include any divisors of log
discrepancy zero into $\Theta$.
\item On the other hand, to be able to apply vanishing, we will need the
supports of $M$ and $\Theta$ to have no components in common.  Thus we cannot allow any divisor
in $\Theta$, which belongs to $M$.
\end{enumerate} 
\end{remark}

Note that $\mathcal{I}_{\Delta,M}$ is a coherent sheaf of ideals of $\ring X.$ which is
determined only by $M$, $\Delta$ and $\mu$.  Since $(K_W+\Theta)-\mu^*(K_X+\Delta)$ and
$\mu ^* \Delta - \Theta$ are effective, one has that
\begin{equation}\label{incl}
\mu _*\ring W.(-\rdown M.)\subset I_{\Delta, M} \subset I_{M} .\tag{$*$}
\end{equation}
We will need the following easy fact:
\begin{lemma}\label{l_sim} Let $\nu \colon \map W'.W.$ be a birational morphism of smooth projective 
varieties, $G$ a reduced divisor on $W$, $\Gamma$ a $\mathbb{Q}$-divisor on $W$ with
$\rdown \Gamma .\leq 0$ and such that $G+\Gamma$ has simple normal crossings support and
where $G,\Gamma$ have no components in common. Then $(W,G+\Gamma )$ is log canonical and
the set of divisors in $W'$ of log discrepancy $0$ for $(W,G+\Gamma )$ and for $(W,G )$
coincide.
\end{lemma}
\begin{proof} See 2.31 of \cite{KM98}.
\end{proof}
\begin{lemma}\label{l_independent} Notation as in \eqref{d_multiplier}.  If $\nu\colon\map W'.W. $ 
is a birational morphism, such that 
$$
\nu ^*M\cup (\nu \circ \mu )^*\Delta \cup \{\text{exceptional locus of }(\nu \circ \mu )\},
$$
has simple normal crossings, then $\mathcal{I}_{\Delta,\nu^*M}$ is also defined and
$\mathcal{I}_{\Delta,M}=\mathcal{I}_{\Delta,\nu^*M}$.
\end{lemma}
\begin{proof} Let $\Theta '$ be the sum of all the divisors in $W'$ of log discrepancy
zero for $(X,\Delta )$. Then $\Theta '$ is also the sum of all divisors in $W'$ of log
discrepancy zero for $(W,\mu ^* \Delta -K_{W/X})$ and so for $(W,\Theta )$ (as one has
$\mu ^* \Delta -K_{W/X}=\Theta -F$ where $F\geq 0$ has no component in common with
$\Theta$ and $\Theta +F$ has simple normal crossings support, and so \eqref{l_sim} applies).
Since $M$ and $\Theta$ have no common components and $M+\Theta$ has simple normal
crossings, one sees that $M$ contains no log canonical centres of $(W,\Theta)$. Thus $\nu
^*M$ has no components in common with $\Theta '$ and so $\mathcal{I}_{\Delta,\nu^*M}$ is
defined.  Note that
\begin{align*} 
M             &= \rdown M.+\{M\} & \text{so that,}\qquad\qquad \\ 
\nu^*M        &= \nu^*\rdown M.+\nu^*\{M\} & \text{and so,}\qquad\qquad\\  
\rdown\nu^*M. &= \nu^*\rdown M.+\rdown \nu^*\{M\}..\\ 
\end{align*} 
By definition 
$$
I_{\Delta,\nu^*M}   = (\nu\circ\mu)_*\ring W'.(K_{W'/X}+\Theta ' -(\nu\circ\mu)^*\Delta-\rdown\nu^*M.)=
$$

$$
\mu_*\left( \ring W.(K_{W/X}+\Theta -\mu ^*\Delta -\rdown M.)\otimes 
\nu _*\ring W'. (K_{W'/W}+\Theta '-\nu ^* \Theta -\rdown \nu^*\{M\}. )\right). 
$$
 Suppose that 
$$
K_{W'/W}+\Theta ' -\nu^*\Theta -\rdown\nu^*\{M\}.,
$$
is effective and exceptional for $\nu$.  Then 
$$
\nu_*\ring W'.(K_{W'/W}+\Theta ' -\nu^*\Theta -\rdown\nu^*\{M\}.)=\ring W.,
$$
and the lemma follows.
We may write
$$
K_{W'}+\Theta ^*=\nu^*(K_{W}+\Theta+\{M\})+E,
$$
where $\Theta ^*$ is the reduced divisor whose components are the divisors on $W'$ of log
discrepancy zero for $(W,\Theta +\{M \} )$ and hence for $(W,\Theta )$. By \eqref{l_sim},
$\Theta ^*= \Theta '$ and so 
$$
K_{W'/W}+\Theta '-\nu ^* \Theta -\rdown \nu ^*\{ M\}. =\rup E.,
$$ 
is effective and exceptional as desired.  \end{proof}

It is important to be able to compare multiplier ideal sheaves.  Here is the first basic
result along these lines:

\begin{lemma}\label{l_compare} Let $(X,\Delta)$ be a smooth log pair, where every 
component of $\Delta$ has coefficient one and let $D$ be a divisor on $X$.  Assume that
there is a birational map $\mu\colon\map W.X.$ with the properties of
\eqref{d_multiplier}.  Let $\Delta'\leq\Delta$ be another divisor, where every component
of $\Delta'$ has coefficient one.

 Then 
$$
\mathcal{I}_{\Delta,M}\subset\mathcal{I}_{\Delta',M}.
$$
\end{lemma}
\begin{proof} Let $F$ be an exceptional divisor extracted by $\mu$.  Then the coefficient of 
$F$ in 
$$
(K_W+\Theta)-\mu^*(K_X+\Delta),
$$
is either the log discrepancy minus one, whenever the log discrepancy is at least one,
or it is zero.  In particular the coefficient of $F$ is an increasing function of the log
discrepancy of $F$.  It follows easily that
$$
(K_W+\Theta)-\mu^*(K_X+\Delta)\leq (K_W+\Theta')-\mu^*(K_X+\Delta'),
$$
which in turn gives the inclusion of ideal sheaves.
\end{proof}

\begin{theorem}\label{t_van} Let $(X,\Delta)$ be a smooth log pair, where every component 
of $\Delta$ has coefficient one and let $\pi\colon\map X.S.$ be a projective morphism of
algebraic varieties.  Let $M$ be an integral divisor such that $M\equiv _{\pi} L+G$ where
\begin{enumerate} 
\item $L$ is $\pi$-nef and there is an effective $\mathbb{Q}$-divisor $B$ such that
$A=L-B$ is $\pi$-ample, 
\item $\Delta$ and $G\cup B$ have no components in common, $\rdown G.=0$ and
\item the support of $\Delta\cup B\cup G$ has simple normal crossings.  
\end{enumerate} 
  
Then $R^p\pi_*\ring X.(K_X+\Delta+M )=0$, for all $p>0$.   
\end{theorem}
\begin{proof} As $L$ is $\pi$-nef, $L-\delta B$ is $\pi$-ample, for any 
$\delta>0$.  
Pick $\epsilon>0$ such that $\epsilon\Delta+L-\delta B$ is $\pi$-ample.  As $\Delta$ and
$G\cup B$ have no components in common and every component of $\Delta$ has coefficient one,
we see that $(X, (1-\epsilon ) \Delta+G+\delta B)$ is kawamata log terminal 
for any $\delta$ sufficiently small.  But then, since
$$\Delta +M\equiv _{\pi} (\epsilon \Delta +L-\delta B)+(1-\epsilon)\Delta +G +\delta B,$$ we may apply Kawamata-Viehweg vanishing.  
\end{proof}

 We fix some notation:

\begin{enumerate} 
\item $(Y,\Gamma)$ is a smooth log pair, 
\item the divisor $X\subset Y$ is a component of $\Gamma$ of coefficient one,
\item $(Y,\Gamma)$ is log canonical, 
\item $S$ is the germ of an algebraic variety, 
\item $\pi\colon\map Y.S.$ is a projective morphism,
\item $\Delta$ is the restriction of $\Gamma-X$ to $X$. 
\end{enumerate} 

Note that we make no requirement about how $X$ sits in $Y$, with respect to $\pi$.  We
will also assume that every coefficient of $\Gamma$ has coefficient one, until \eqref{t_asymptotic}.  

\begin{definition}\label{d_pi} A divisor $D$ on $Y$ will be called 
\textbf{$\pi$-effective for the pair} 
$(Y, X)$ if the natural homomorphism
$$
\map {\pi_*\ring Y.(D)}.{\pi_*\ring X.(D)}.,
$$
is not zero.  It is called \textbf{$\pi$-$\mathbb{Q}$-effective for the pair} $(Y,X)$
if $mD$ is $\pi$-effective for the pair $(Y, X)$ for some positive integer $m$.  $D$ is
said to be \textbf{$\pi$-big for the pair} $(Y,X)$ if we can write $mD=A+B$ for a positive
integer $m$, a $\pi$-ample divisor $A$ and a $\pi$-effective divisor $B$ for the pair
$(Y,X)$.

A divisor $D$ on $X$ will be called \textbf{$\pi$-transverse} for $(X,\Delta )$ if the
sheaf $\ring X.(D)$ is $\pi$-generated at the generic point of every log canonical centre
of $K_X+\Delta$.  $D$ will be called \textbf{$\pi$-$\mathbb{Q}$-transverse} for $(X,\Delta
)$ if $mD$ is $\pi$-transverse for some $m>0$.
\end{definition}

Now we turn to the definition of sheaves $\mathcal{J}^0_{\Delta,D}$,
$\mathcal{I}^0_{\Delta,D}$, $\mathcal{J}^1_{\Delta,D}$, $\mathcal{I}^1_{\Delta,D}$, which
generalize the sheaves $\mathcal{J}^0_D$, $\mathcal{I}^0_D$,
$\mathcal{J}^1_D$, $\mathcal{I}^1_D$ of Definition 2.10 of \cite{Kawamata99}.   

 We will need a:
\begin{definition}\label{d_canonical} Let $(X,\Delta)$ be a log pair and 
let $\mu\colon\map W.X.$ be a birational morphism.  We say that $\mu$ is
\textbf{canonical} if every exceptional divisor extracted by $\mu$ has log discrepancy at
least one.
\end{definition}

If $(X,\Delta)$ is a pair consisting of a smooth variety and a divisor with simple normal
crossings support, then note that $\mu$ is canonical provided it is an isomorphism in a
neighbourhood of the generic point of any irreducible component of the intersection of the
components of $\Delta$, that is provided it is an isomorphism in a neighbourhood of any log
canonical centre.  Note also that if $\mu$ is canonical, then $\Theta$ is the strict
transform of $\Delta$.

\begin{definition}\label{d_ideal} Let $D$ be a $\pi$-$\mathbb{Q}$-transverse divisor for 
$(X,\Delta )$.  For each positive integer $m$, such that $mD$ is $\pi$-transverse, pick a
proper birational canonical morphism $\mu_m\colon\map W_m.X.$ from a smooth variety $W_m$,
such that
\begin{enumerate} 
\item there is a decomposition $\mu^*_{m}(mD)=P_m+M_m$ in $\Div(W_m)$ such that the natural homomorphism
$$
\map {(\pi\circ\mu_m)_*\ring W_m.(P_m)}.{\pi_*\ring X.(mD)}.,
$$
is an isomorphism,  
\item $P_m$ is $\pi\circ\mu_m$-free, and
\item $M_m\cup \Theta_m \cup\{\text{the exceptional set of }\mu _m\}$ is effective, with normal crossing support,
where $\Theta_m$ is defined as in \eqref{d_multiplier}.  
\end{enumerate} 
We define
$$
\mathcal{J}^0_{\Delta,D}=\bigcup \mathcal{I}_{\Delta,\frac 1mM_m},
$$
where the union is taken over all positive integers $m$ such that $mD$ is
$\pi$-transverse.  If $D$ itself is $\pi$-transverse then we set
$\mathcal{I}^0_{\Delta,D}=\mathcal{I}_{\Delta,M_1}$.

Now suppose that $D$ is a $\pi$-$\mathbb{Q}$-transverse divisor on $Y$ for $(Y,\Gamma
)$.  For every positive integer $m$, such that $mD$ is a $\pi$-transverse divisor for
$(Y,\Gamma )$, pick a birational morphism $\mu_m\colon\map W_m.Y.$ from a smooth variety
$W_m$, such that
\begin{enumerate} 
\item there is a decomposition $\mu^*_m(mD)=Q_m+N_m$ in $\Div(W_m)$ such that the natural homomorphism
$$
\map {(\pi\circ\mu_m)_*\ring W_m.(Q_m)}.{\pi_*\ring Y.(mD)}.,
$$
is an isomorphism,  
\item $Q_m$ is $\pi\circ\mu_m$-free, and
\item $N_m$ is effective, and $X_m+N_m+\Theta_m +\{\text{the exceptional set of }\mu _m\}$
has normal crossing support, where $X_m$ is the strict transform of $X$ and $\Theta_m$
is defined as in \eqref{d_multiplier}.
\end{enumerate} 
By assumption $X_m$ is not contained in the support of $N_m$.  We define
$$
\mathcal{J}^1_{\Delta,D}=\bigcup \mathcal{I}_{\Delta,\frac 1mN_m|_{X_m}},
$$
where the union is taken over all positive integers $m$ such that $mD$ is
$\pi$-transverse for the pair $(Y,\Gamma )$.  If $D$ itself is $\pi$-transverse for 
$(Y,\Gamma)$, then we set $\mathcal{I}^1_{\Delta,D}=\mathcal{I}_{\Delta,N_1|_{X_1}}$.
\end{definition}

\begin{remark}\label{r_notneed} We will not make use of the sheaves $\mathcal{I}^i_{\Delta,D}$, 
for $i=0$ and $1$.  We have included their definition for the sake of completeness.
\end{remark}

\begin{lemma}\label{l_can} Let $D$ be a $\pi$-$\mathbb{Q}$-transverse 
divisor for $(Y,\Gamma )$ (respectively for $(X,\Delta )$). Then the sheaf $\mathcal{J}^i
_{\Delta,D}$ is defined, it does not depend on the resolution $\mu _m$ and we may choose
$\mu _m$ to be canonical for any $m$ sufficiently divisible, where $i=1$ (respectively
$i=0$).  Further, we may assume that there is a $\mu _m$-ample divisor whose support is
exceptional.
\end{lemma}
\begin{proof} 
This is an easy consequence of \eqref{l_independent} and of the Resolution Lemma of \cite{Szabo94}.
\end{proof}

\begin{lemma}\label{l_inclusion} Let $D$ be a $\pi$-$\mathbb{Q}$-transverse 
divisor for $(Y,\Gamma )$.  
\begin{enumerate} 
\item $\mathcal{J}^1_{\Delta,D}\subset \mathcal{J}^0_{\Delta,D|_X}$
\item If $B$ is an effective integral divisor, which does not contain $X$
and $D+B$ is $\pi$-$\mathbb{Q}$-transverse for $(Y,\Gamma)$, then 
$$
\mathcal{J}^i_{\Delta,D}(-B)\subset \mathcal{J}^i_{\Delta,D+B}, \qquad \text{$i=0$, $1$.}
$$
\item There is a positive integer $m$ such that
$$
\mathcal{J}^0_{\Delta,D}=\mathcal{I}_{\Delta,\frac 1mM_m} \qquad \text{and} \qquad \mathcal{J}^1_{\Delta,D}=\mathcal{I}_{\Delta,\frac 1mN_m|_{X_m}}.
$$
\end{enumerate} 
\end{lemma}
\begin{proof} (1) follows, since $N_m|_{X_m}\geq M_m$ as divisors on $X_m$.  (2) is clear.  (3) follows as $X$
is Noetherian.
\end{proof}

\begin{lemma}\label{l_basic} Let $D$ be a $\pi$-$\mathbb{Q}$-transverse divisor for  
$(Y , X, \Gamma)$.  
Then 
\begin{enumerate} 
\item $\mathcal{J}^i_{\Delta,\alpha D}\subset \mathcal{J}^i_{\Delta,D}$, for any rational number $\alpha\geq 1$ and 
$i=0$, $1$.     
\item $\mathcal{J}^i_{\Delta,D}\subset \mathcal{J}^i_{\Delta,D+L}$, for any $\pi$-free divisor $L$, and $i=0$, $1$.        
\item 
$
\im(\map\pi_*{\ring Y.(D)}.{\pi_*\ring X.(D)}.)\subset \pi_*\mathcal{J}^1_{\Delta,D}(D).
$
\end{enumerate} 
\end{lemma}
\begin{proof} First of all, notice that Definition \eqref{d_ideal} also makes sense for $\mathbb{Q}$-divisors. (1) and (2) are clear.  We may assume that 
$\mu _m$ is chosen so that $\mu _m ^*D=P+M$ in $\Div(W_m)$ and
$$
(\pi \circ \mu _m)_*\map {\ring W_m.}(P).{\pi _*\ring Y.(D)}.,
$$ 
is an isomorphism.  Since $M>M_m/m$ one sees that there are inclusions
$$
\ring W_m.(\mu_m^*D-M)\subset \ring W_m.(\mu _m^* D-\rdown M_m/m.)\subset \ring W_m.(\mu _m^* D ).
$$ 

We now push forward via $(\pi\circ\mu _m)$.  Since $P=\mu_m^* D-M$, both the left and
right hand sides push forward to $\pi _* \ring Y.(D)$.  So the image of $\pi _* \ring
Y.(D)$ is contained in the image of $(\pi \circ \mu _m)_*\ring W_m.(\mu _m^* D-\rdown
M_m/m.)$ and hence in $\pi_*\mathcal{J}^1_{\Delta,D}(D)$ (cf. equation \eqref{incl}, after
\eqref{r_multiplier}).
\end{proof}
The following lemma is the key result that will allow us to generalize the
results of Siu and Kawamata to the setting of log-forms:
\begin{lemma}\label{l_include} Suppose that $D$ is a divisor on $Y$ and suppose that 
there is an effective divisor $B$ on $Y$ such that 
\begin{enumerate} 
\item $A=D-B$ is $\pi$-ample, and
\item $B$ is $\pi $-$\mathbb{Q}$-transverse for $(Y,\Gamma )$.
\end{enumerate} 
 Then 
$$
\pi_*\mathcal{J}^1_{\Delta,D}(D+K_X+\Delta)\subset \im(\map\pi_*{\ring Y.(D+K_Y+\Gamma)}.\pi_*{\ring X.(D+K_X+\Delta)}.).
$$
\end{lemma}
\begin{proof} Notation as in \eqref{d_ideal}.  There is an exact sequence
\begin{multline*}
0\longrightarrow \ring W_m.(K_{W_m}+\Theta_m-X_m+\mu _m^*D-\rdown\frac 
{N_m}{m}.)\longrightarrow \ring W_m.(K_{W_m}+\Theta_m+\mu _m^*D-\rdown\frac 
{N_m}{m}
.) \\
\longrightarrow\ring X_m.(K_{X_m}+(\Theta_m-X_m)|_{X_m}+(\mu _m^*D-\rdown\frac 
{N_m}{m}.)|_{X_m})\longrightarrow 0. 
\end{multline*}
\begin{claim} If $m$ is sufficiently large then
$$
R^1(\pi\circ\mu _m)_*(\ring W_m.(K_{W_m}+\Theta_m-X_m+\mu _m^*D-\rdown\frac {N_m}{m}.))=0.
$$
\end{claim}
\begin{proof}[Proof of Claim] By \eqref{l_can}, we may assume that $\mu _m$ is
canonical and there is an effective divisor $E$ supported on the exceptional locus such
that $\mu _m^*A-E$ is $\pi \circ \mu _m$-ample. As $\mu _m$ is canonical, $E+\{ N_m/m\}$
and $\Theta _m -X_m$ have no common component. Since $B$ is $\pi$-$\mathbb{Q}$-transverse, we may
assume that $\mu _m ^* (mB)\sim _{\pi, \mathbb{Q}}N'$ where $N'$ and $\Theta _m -X_m$ have no
common component and $N_m+N'+\Theta _m+\{\text{exceptional locus of } \mu _m\}$ has simple
normal crossings support.  Notice that since $D-B$ is $\pi$-ample, we have that $D$ is
$\pi$-$\mathbb{Q}$-transverse and that $N_m\leq N'$.  We now apply \eqref{t_van} to $\Delta
=\Theta _m-X_m$, $M=\mu _m^*D- \rdown N_m/m.$, $G=\{ N_m/m\}$ and $L=Q_m/m$. Notice that
$L-(\frac{N'-N_m}{m}+E)\sim _{\pi \circ \mu _m}\mu _m^*A -E$ is $\pi\circ \mu _m$-ample
and $\frac{N'-N_m}{m}+E$ is effective.  The hypothesis of \eqref{t_van} are now verified
and so
$$R^p(\pi \circ \mu _m)_*\ring W_m.(K_{W_m}+\Theta_m -X_m +\mu _m^*D-\rdown\frac {N_m}{m}.)$$
vanishes for $p>0$.
\end{proof}
Thus there is a surjective homomorphism
\begin{multline*}
\pi_*\ring Y.(D+K_Y+\Gamma)\supset (\pi\circ\mu _m)_*\ring W_m.(K_{W_m}+\Theta_m+\mu _m^*D-\rdown\frac {N_m}{m}.)\\
\map .{(\pi\circ\mu _m)_*\ring X_m.(K_{X_m}+\Theta_m-X_m+\mu _m^*D-\rdown\frac {N_m}{m}.)}.=\mathcal{J}^1_{\Delta,D}(D+K_X+\Delta),
\end{multline*}
hence the assertion.  
\end{proof}

\begin{theorem}\label{t_asymptotic} Let $X\subset Y$ be a smooth divisor in a smooth 
variety. Fix a positive integer $e$.  Let $H$ be a sufficiently $\pi$-very ample divisor
and set $A=(\dim X+1)H$.  Assume that
\begin{enumerate}
\item $\Gamma$ is a $\mathbb{Q}$-divisor with simple normal crossings support such that
$\Gamma$ contains $X$ with coefficient one and $(Y,\Gamma )$ is log canonical;
\item $k$ is a positive  integer such that $k\Gamma$ is integral;
\item $C\geq 0$ is an integral divisor not containing $X$;
\item Given $D=k(K_Y+\Gamma)$, $D'=D|_X=k(K_X+\Delta)$ and $G=D+C$, then $D'$ is
$\pi$-$\mathbb{Q}$-pseudo-effective (i.e. $D'+\epsilon H'$ is
$\pi$-$\mathbb{Q}$-effective for any $\epsilon >0$) and that $G$ is
$\pi$-$\mathbb{Q}$-transverse for $(Y,\rup \Gamma.)$.
\end{enumerate}
Then
$$
\mathcal{J}^0_{mD'+eH'}(-mC)\subset \mathcal{J}^1_{mG+eH+A}, \qquad \forall m>0,
$$
where $H':=H|_X$.
\end{theorem}
\begin{proof} It suffices to prove the stronger inclusion
$$
\mathcal{J}^0_{mD'+eH'}(-mC)\subset \mathcal{J}^1_{\rup\Delta.,mG+eH+A}.
$$
This is a modification of the proof of Lemma 3.7 of \cite{Kawamata99}, which we follow
quite closely.  We prove the inclusion by induction on $m$, in the case $e=1$, the general
case being left to the reader.

By assumption $mD'+H'$ is $\pi$-big, $C$ does not contain $X$ and $mG+H+A$ is
$\pi$-$\mathbb{Q}$-transverse for $(Y,\rup \Gamma . )$ for any non-negative integer $m$.
In particular, both sides of the proposed inclusion are well-defined.

If $m=0$ the result is clear.  Assume the result for $m$.  Suppose that 
$$
\Gamma=\sum _j \gamma_j \Gamma(j), \qquad \text{where} \qquad 0< \gamma_j :=\frac {g_j}k\leq 1.
$$
Then set 
$$
\Gamma_i=\sum _j \gamma_i(j)\Gamma(j)\qquad \text{where} \qquad \gamma_i(j)=\begin{cases} 0 & \text{if $1\leq i\leq k-g_j$,} \\
                                                                                          1 & \text{if $k-g_j<i\leq k+1$.} \end{cases}
$$
With this choice of $\Gamma_i$, we have
$$
X\leq \Gamma_1\leq \Gamma_2\leq \Gamma_3\leq \dots \leq \Gamma_k=\Gamma_{k+1}=\rup\Gamma.,
$$
and each $\Gamma _i$ is integral.  Let $\Delta_i$ be the restriction of $\Gamma_i-X$ to
$X$.  Set $$D_i=K_Y+\Gamma_i\qquad {\text {and}} \qquad D_{\leq i }=\sum _{j\leq i}D_j.$$
With this choice of $D_i$, we have $D=D_{\leq k}$.  Possibly replacing $H$ by a multiple,
we may assume that $H_i=D_{\leq i-1}+\Gamma _i-X+H$ and $A _i= H+D_{\leq i}+A$ are
$\pi$-ample, for $1\leq i\leq k$.  We are now going to prove, by induction on $i$, that
$$
\mathcal{J}^0_{mD'+H'}(-mC)\subset \mathcal{J}^1_{\Delta_{i+1},mG+D_{\leq i}+H+A} \qquad \text{where $0\leq i\leq k$,}
$$
where we adopt the convention that $D_{\leq 0}=0$.  

 Since $\Delta_1\leq\rup\Delta.$, it follows by \eqref{l_compare}
$$
\mathcal{J}^1_{\rup\Delta.,mG+H+A}\subset\mathcal{J}^1_{\Delta_1,mG+H+A},
$$
and so the case $i=0$ follows from the inclusion
$$
\mathcal{J}^0_{mD'+H'}(-mC)\subset \mathcal{J}^1_{\rup\Delta.,mG+H+A},
$$
which we are assuming by induction on $m$.  Assume the result up to $i-1$.  Then
\begin{align*} 
&\pi_*\mathcal{J}^0_{mD'+H'}(mD+D_{\leq i}+H+A) \\
&=\pi_*\mathcal{J}^0_{mD'+H'}(-mC)(mG+D_{\leq i}+H+A) \\
&\subset\pi_*\mathcal{J}^1_{\Delta_i,mG+D_{\leq i-1}+H+A}(mG+D_{\leq i}+H+A) \\
&=\pi_*\mathcal{J}^1_{\Delta_i,mG+D_{\leq i-1}+H+A}(mG+D_{\leq i-1}+H+A+K_Y+\Gamma_i) \\
&\subset\im(\map {\pi_*\ring Y.(mG+D_{\leq i}+H+A)}.{\pi_*\ring X.(mG+D_{\leq i}+H+A)}.)\\
&\subset \pi_*\mathcal{J}^1_{\Delta_{i+1},mG+D_{\leq i}+H+A}(mG+D_{\leq i}+H+A), \\
\end{align*}
where we use induction to get the inclusion of the second line in the third line,
\eqref{l_include} applied to the basic identity
\begin{align*} 
mG+D_{\leq i-1}+H+A  &=  A_{i-1}+B,\\ 
\end{align*} 
to get the inclusion of the fourth line in the fifth line, and (3) of \eqref{l_basic}
to get the inclusion of the fifth line in the sixth line.  As
$$
\mathcal{J}^0_{mD'+H'}(mD+D_{\leq i}+H+A)=\mathcal{J}^0_{mD'+H'}((mD+H)+H_i+(\dim X)H+K_X),
$$
is $\pi$-generated by (a trivial generalization of) (3.2) of \cite{Kawamata99}, it follows that
$$
\mathcal{J}^0_{mD'+H'}(-mC)\subset \mathcal {J}^1_{\Delta_{i+1},mG+D_{\leq i}+H+A}.
$$
This completes the induction on $i$. It follows that 
$$
\mathcal{J}^0_{mD'+H'}(-mC)\subset \mathcal{J}^1_{\Delta_{k+1},mG+D_{\leq k}+H+A},
$$ 
so that, by (2) of \eqref{l_inclusion}, 
$$
\mathcal{J}^0_{mD'+H'}(-(m+1)C)\subset \mathcal{J}^1_{\rup\Delta.,(m+1)G+H+A}. 
$$ 
But by (1) and (2) of \eqref{l_basic}, 
$$
\mathcal{J}^0_{(m+1)D'+H'}\subset \mathcal{J}^0_{mD'+H'}.
$$
Thus 
$$
\mathcal{J}^0_{(m+1)D'+H'}(-(m+1)C)\subset \mathcal{J}^1_{\rup\Delta.,(m+1)G+H+A}, 
$$
and this completes the induction on $m$ and the proof.  
\end{proof}

\begin{corollary}\label{c_ext} Let $X\subset Y$ be a smooth divisor in a smooth variety.  
Let $H$ be a sufficiently $\pi$-very ample divisor and set $A=(\dim X+1)H$.  Assume that
\begin{enumerate}
\item $\Gamma$ is a $\mathbb{Q}$-divisor with simple normal crossings support such that
$\Gamma$ contains $X$ with coefficient one and $(Y,\Gamma )$ is log canonical;
\item $C\geq 0$ is a $\mathbb{Q}$-divisor not containing $X$;
\item $K_X+\Delta$ is $\pi$-$\mathbb{Q}$-pseudo-effective, where $\Delta=(\Gamma-X)|_X$
and $K_Y+\Gamma+C$ is $\pi$-$\mathbb{Q}$-transverse for $(Y,\rup \Gamma.)$.
\end{enumerate}
For any positive integer $m$, such that $m(K_Y+\Gamma)$ is integral, the image of the
natural homomorphism
$$
\map {\pi _* \ring Y.(m(K_Y+\Gamma+C)+H+A)}.{ \pi _* \ring X.(m(K_X+\Delta+C)+H+A)}.,
$$
contains the image of the sheaf $\pi _*\ring X.(m(K_X+\Delta)+H)$ considered as a subsheaf of $\pi _*
\ring X.(m(K_X+\Delta+C)+H+A)$ by the inclusion induced by any divisor in $mC+|A|$ not containing $X$.
\end{corollary}
\begin{proof} Let $H'=H|_X$.  Fixing a section $\sigma \in H^0(Y,\ring Y.(A))$ not vanishing on $X$, 
we can view 
$$
\mathcal{J} ^0_{m(K_X+\Delta)+H'}(m(K_Y+\Gamma)+H),
$$
as a subsheaf of
$$
\mathcal{J}^0_{m(K_X+\Delta)+H'}(m(K_Y+\Gamma)+H+A),
$$
via the map induced by multiplication by $\sigma$.  As we have seen in the proof of 
\eqref{t_asymptotic}, 
\begin{multline*}
\pi _*\mathcal{J}^0_{m(K_X+\Delta)+H'}(m(K_Y+\Gamma)+H+A)= \\
\pi _*\mathcal{J}^0_{m(K_X+\Delta)+H'}(-mC)(m(K_Y+\Gamma+C)+H+A),
\end{multline*}
is contained in the image of $\pi _* \ring Y.(m(K_Y+\Gamma+C)+H+A)$.  Since 
$$
\pi _* \ring X.(m(K_X+\Delta)+H)=\pi _*\mathcal{J}^0_{m(K_X+\Delta)+H'}(m(K_Y+\Gamma)+H),
$$
the assertion now follows easily.
\end{proof}

\section{Lifting log canonical centres}
\label{s_centres}

We fix some notation for this section.  Let $(X,\Delta)$ be a log pair, where $X$ is
$\mathbb{Q}$-factorial, and $\Delta$ is an effective $\mathbb{Q}$-divisor.  Let $V$ be an
exceptional log canonical centre of $K_X+\Delta$.  Let $f\colon\map W.V.$ be a resolution
of $V$.  Let $\Theta$ be an effective $\mathbb{Q}$-divisor on $W$.  Suppose that there are
positive rational numbers $\lambda$ and $\mu$ such that $\Delta\sim\lambda K_X$ and
$\Theta\sim\mu K_W$.  Let $\nu=(\lambda+1)(\mu+1)-1$.  Suppose that $W$ is of general
type.

The main result of this section is:

\begin{theorem}\label{t_lift} There is a very general subset $U$ of $V$ with the following
property:

Suppose that $W'\subset W$ is a pure log canonical centre of $K_W+\Theta$, whose image
$V'\subset V$ intersects $U$.

Then, for every positive rational number $\delta$, we may find an effective divisor
$\Delta'$ on $X$ such that $V'$ is a pure log canonical centre of $K_X+\Delta'$, where
$\Delta'\sim(\nu+\delta)K_X$.  Now suppose that we may write $K_X\sim A+E$, where $A$ is
ample, $E$ is effective and $V$ is not contained in $E$.  Then we may assume that $V'$ is an
exceptional log canonical centre of $K_X+\Delta'$. 
\end{theorem}

\eqref{t_lift} is the main step in the proof of \eqref{t_birational}.  In fact
\eqref{t_birational} is a standard consequences of \eqref{t_lift}, see \S \ref{s_proof}.

The idea is to lift $\Theta$ as a $\mathbb{Q}$-divisor, using \eqref{c_ext}.  However this
is more delicate than it might first appear, as we cannot lift $\Theta$ directly to $X$.
Instead, we are able to find a sequence of successively better approximations $\Theta_m$.
In particular we have to pass to a resolution of $X$, and to this end, we introduce some
more notation.

Let $\pi\colon\map Y.X.$ be a log resolution of $(X,\Delta)$.  Let $E$ be the unique
exceptional divisor of log discrepancy zero with centre $V$.  Possibly blowing up further,
we may assume that the restriction $p\colon\map E.V.$ of $\pi$ to $E$ factors $g\colon\map
E.W.$ through $f$.  Thus we have a commutative diagram,
$$
\begin{diagram}
E_{\xi}& \subset  &  E  &   \rInto    & Y   \\
&\ldTo^{g}   &     &           &      \\
W&  & \dTo^p    &           & \dTo_{\pi} \\
&\rdTo_{f}  &     &           &      \\
\xi & \in  &  V  &   \rInto  & X.
\end{diagram}
$$ 
Let $\xi$ be the generic point of $V$ and let $E_{\xi}$ be the generic fibre of $p$.  

\begin{lemma}\label{l_limiting} We may find a positive integer $m_0$ and a positive integer
$s$, with the following properties:

For every sufficiently large and divisible integer $m$, there is an effective divisor
$G_m$ on $Y$ such that
\begin{enumerate} 
\item $G_m$ does not contain $E$,
\item $G_m+sE\sim (m_0+m(1+\lambda))\pi^*K_X$, 
\item $G_m|_E=mg^*\Theta+B_m$, where the restriction of $B_m$ to the inverse image of a
fixed open subset of $V$ is effective, and
\item $B_m|_{E_{\xi}}$ belongs to a fixed linear system.
\end{enumerate} 
\end{lemma}
\begin{proof} We may assume that 
$$
K_Y+\Gamma =\pi ^*(K_X+\Delta)+F,
$$
where $\Gamma$ and $F$ are effective, with no common components and $F$ is exceptional.  We
may write $\Gamma=E+\Gamma ^h+\Gamma ^v$, $F=F^h+F^v$ where the components of $\Gamma ^h$
and $F^h$ are the components of $\Gamma -E$ and $F$ whose images contain $V$.  As $V$ is
an exceptional log canonical centre for $(X,\Delta)$, $\{\Gamma^h\}=\Gamma ^h$ and
$\Gamma ^h |_{E_{\xi}}= (\Gamma -E )|_{E_{\xi}}$ and so, 
$$
(K_E+\Gamma ^h)|_{E_{\xi}}=(K_Y+\Gamma)|_{E_{\xi}} = F|_{E_{\xi}}\geq  0.
$$

By Corollary \eqref{c_wp}, we may assume that 
$$
h^0(E,\ring E.(m(K_{E/W}+\Gamma^h)+H))>0,
$$
where $H$ is a sufficiently ample divisor on $Y$.  Therefore, there is an injection
$$
\map g^*|mK_W|.|m(K_E+\Gamma ^h)+H|.,
$$
induced by a choice of a divisor in $|m(K_{E/W}+\Gamma^h)+H|$.  

We now apply \eqref{c_ext}.  Possibly blowing up further, we may assume that the
image of every element of $\lcc(Y,E+\rup \Gamma ^h.)$ contains $V$.  As the divisors
$E+\rup \Gamma ^h.$ and $F$ have no common components, it follows that
$K_Y+\Gamma=\pi^*(K_X+\Delta)+F$ is $\pi$-$\mathbb{Q}$-transverse for $(Y,E+\rup \Gamma^h.)$.
Notice also that as $W$ is of general type, $(K_Y+E+\Gamma^h)|_E=(K_E+\Gamma^h|_E)$ is
$\mathbb{Q}$-effective, by \eqref{c_nm}.  By Corollary \eqref{c_ext}, with
$S=\sp(\mathbb{C})$, it follows that elements of
$$
\emb |m(K_E+\Gamma^h|_E)+H|_E|.|m(K_E+(\Gamma-E)|_E)+H|_E+A|_E|.,
$$
can be lifted to elements of 
$$
|m(K_Y+\Gamma)+H+A|=|m(\pi ^*(K_X+\Delta)+F)+H+A|.
$$

Denote by $S_m$ the lift of $mg^*\Theta$ to $|m\mu(\pi ^*(K_X+\Delta)+F)+H+A|$.  As $K_X$ is
big, there is a fixed positive integer $m_0$ and an effective divisor
$$
S_0\in |m_0\pi^* K_X-H-A|. 
$$
Then 
$$
S_m+S_0\in |(m_0+m(1+\lambda)\mu)\pi ^*K_X+m\mu F|=|(m_0+m(1+\lambda)\mu)\pi^*K_X|+m\mu F
$$
so that we may write $S_m+S_0=G'_m+m\mu F$ for some 
$$
G'_m\in |(m_0+m(1+\lambda)\mu)\pi^*K_X|.
$$
We may write 
$$
G'_m=G_m+sE,
$$
where $G_m$ does not contain $E$ and $s$ is a positive integer.  Since $S_m$ does not
contain $E$ and $S_0$ is a fixed divisor, it is clear that $s$ is fixed.  Since $G_m+m\mu F$
contains $S_m$, it follows that we may write
$$
G_m|_E=mg^*\Theta+B_m,
$$
where the only components of $B_m$ with negative coefficients, are supported on
$g^*\Theta$ and contained in $E\cap F$.  Thus (3) holds, and when we restrict to $E_{\xi}$
we may ignore $g^*\Theta$.  Thus
$$
B_m|_{E_{\xi}}\sim ((m_0+m(1+\lambda)\mu)\pi^*K_X-sE)|_{E_{\xi}}=(-sE)|_{E_{\xi}}.
$$
\end{proof}

\begin{proof}[Proof of \eqref{t_lift}] 

Let $U$ be any very general subset of $V$, such
that
\begin{enumerate} 
\item the morphism $p$ restricted to $p^{-1}(U)$ is smooth, 
\item every exceptional divisor whose centre intersects $U$, contains $V$, 
\item every component of $\Delta$ which intersects $U$, contains $V$, and
\item every component of $B_m$ whose image intersects $U$, dominates $V$, 
for every sufficiently large and divisible integer $m$, as in \eqref{l_limiting}.
\end{enumerate} 

Let $G_m$ be the sequence of divisors whose existence is guaranteed by \eqref{l_limiting}.
Set
$$
\Gamma_m=G'_m/m, \qquad \Theta_m=\pi_*\Gamma_m \qquad \text{and} \qquad \Delta_m=a_m\Delta+\Theta_m,
$$
where $a_m\leq 1$ is chosen so that the log discrepancy of $E$ with respect to
$K_X+\Delta_m$ is zero.

\begin{claim}\label{c_irr} $V'$ is a pure log canonical centre of $K_X+\Delta_m$.
\end{claim}
\begin{proof}[Proof of Claim] Since this result is local about the generic point of $V$,
passing to an open subset of $X$, we may assume that properties (1-4) hold, and that $B_m$
is effective.  In this case
\begin{align*} 
K_Y+E+F_m+G_m/m &=\pi^*(K_X+\Delta_m)\qquad \text{where}  \\ 
K_Y+c_mE+F_m    &=\pi^*(K_X+a_m\Delta), \\ 
\end{align*} 
$F_m$ is a divisor, which does not contain $E$ and $c_m=1-s/m<1$ is a positive number.
Restricting to $E$, we get
$$
K_E+F_m|_E+g^*\Theta+B_m/m.
$$
Let $E'$ be the inverse image of $W'$.  Since the morphism $p$ is smooth, $E'$ is a
pure log canonical centre of $K_E+g^*\Theta$.  Since $V$ is an exceptional log canonical
centre, we have $\rdown \Gamma^h.=0$.  Thus $K_E+\Gamma^h$ is kawamata log terminal.  Since
$B_m|_{E_{\xi}}$ belongs to a fixed linear system and $p$ is smooth over $U$, it follows
that the restriction of $B_m$ to any fibre of $p$ over $U$, belongs to a fixed linear
system.  Thus $K_E+\Gamma^h|_E+B_m/m$ is kawamata log terminal, for $m$ large enough.  As
$F_m\leq \Gamma^h$, it follows that $K_E+F_m|_E+B_m/m$ is also kawamata log
terminal, for $m$ large enough.  It follows that $E'$ is a pure log canonical centre of
$K_E+F_m|_E+g^*\Theta+B_m/m$.

We now apply inversion of adjunction, see (17.7) of \cite{Kollaretal}.  (17.7) is only
stated for effective divisors.  However it is easy to see that (17.6) and (17.7) of
\cite{Kollaretal}, extend to the case when the components of negative coefficient are
contracted by a birational map, since this is the case for (17.4) of \cite{Kollaretal}.
\end{proof}

Clearly $0\leq a_m\leq 1$.  Moreover, $a_m$ approaches $1$ as $m$ approaches infinity, as
the coefficient $s/m$ of $E$ in $\Gamma_m$ approaches zero.  Pick $m$ so that
$m_0/m<\delta/2$.  Then, using the decomposition $K_X\sim A+E$, there is a perturbation
$\Delta'\sim(\nu+\delta)K_X$ of $\Delta_m$, such that $V'$ is an exceptional log canonical
centre of $K_X+\Delta'$.
\end{proof}

\section{Effective birational freeness}
\label{s_effective}

In this section we give an effective bound for the $r$th pluricanonical map of a variety
of general type, which only depends on the volume of $K_X$ and the value of this constant
for general subvarieties of $X$:

\begin{proposition}\label{p_effective} Let $X$ be a smooth variety of general type 
of dimension $n$.  Let $s'$ be a positive integer such that for every subvariety $V$ of
$X$ which contains a very general point of $X$, the map $\phi_{s'K_V}$ is birational.  Set
$$
s=\left (\rdown \frac{n}{\vol (K_X)^{1/n}}.+2\right )\left (ns'+1\right ).
$$

 Then $\phi_{rK_X}$ is birational, for all $r\geq 4s(s'+1)+s'$.  
\end{proposition}
\begin{proof} By \eqref{l_sing} and \eqref{l_d}, for every very general point $x\in X$, we 
may find a $\mathbb{Q}$-divisor $\Delta\sim _{\mathbb{Q}}\lambda K_X$ with
$$
\lambda < \rdown \frac n{\vol (K_X)^{1/n}}.+1,
$$
such that $\lcc(X,\Delta,x)\ne \emptyset$.  By \eqref{l_irr}, we may assume that
$(X,\Delta)$ is log canonical near $x$, that $\lcc (X,\Delta,x)$ contains a unique element
$V$ and that there is a unique place of log canonical singularities dominating $V$.  Let
$f\colon\map W.V.$ be a resolution of $V$.  Then, by \eqref{l_birational}, for a general
point $x'\in V$, we may find a divisor $\Theta$ on $W$, such that $x'$ is an isolated log
canonical centre of $K_W+\Theta$, where $\Theta\sim_{\mathbb{Q}} ns' K_W$.  By
\eqref{t_lift} it follows that we may find $\Delta'\sim (s-1)K_X$ such that $x$ is an
isolated point of $K_X+\Delta'$.  By \eqref{l_lift} it follows that $h^0(2sK_X)\geq 2$,
and so the result follows from \eqref{l_bir}.
\end{proof}

\begin{corollary}\label{c_lower} Let $X$ be a smooth variety of general type 
of dimension $n$.  Let $s'$ be a positive integer such that for every subvariety $V$ of
$X$ which contains a very general point of $X$, the map $\phi_{s'K_V}$ is birational.

 If the volume of $K_X$ is at most $v$, then $X$ is birational to a subvariety of 
degree at most
$$
(8(ns'+1)(s'+1))^n \max \left (n^n,3^n v \right ).
$$
\end{corollary}
\begin{proof} By \eqref{p_effective}, $\phi_{rK_X}$ is birational, 
where 
$$
r=4s(s'+1)+s' \qquad \text{and} \qquad s=\left (\rdown \frac{n}{\vol (K_X)^{1/n}}.+2\right )\left (ns'+1\right ).
$$
Now 
\begin{align*} 
\vol(rK_X) &= r^n\vol(K_X) \\ 
           &=(4s(s'+1)+s')^n\vol(K_X) \\
           &\leq ((4s+1)(s'+1))^n\vol(K_X) \\
           &\leq \left (\rdown \frac{n}{\vol (K_X)^{1/n}}.+3\right )^n (4(ns'+1)(s'+1))^n\vol(K_X)  \\ 
           &\leq \max \left (\frac {n^n}{\vol(K_X)},3^n\right ) (8(ns'+1)(s'+1))^n\vol(K_X)  \\ 
           &\leq (8(ns'+1)(s'+1))^n \max \left (n^n,3^n v \right ), \\ 
\end{align*} 
and so we are done by  \eqref{l_low}.  
\end{proof}

\section{Proof of \eqref{t_birational} }
\label{s_proof}

\begin{lemma}\label{l_bb} Let $\pi \colon\map X.B.$ be a bounded family of projective 
varieties of general type. Then there is a positive integer $R$ such
that if $Y$ is a resolution of any fiber of $\pi$, then $\phi _{rK_Y}$
is birational for all $r\geq R$.
\end{lemma}
\begin{proof} Let $Y$ be a fiber over the generic point $\xi$ of a component of
$B$. If we pick a projective resolution of $Y$, then we may extend this to an open
neighbourhood of $\xi$. Thus, by Noetherian induction, possibly replacing $B$ by a union of
locally closed subsets, we may assume that $\pi$ is projective and smooth. If $\phi
_{rK_Y}$ is birational, then this map remains birational over an open neighbourhood of
$\xi$, and we are done by Noetherian induction.
\end{proof}

\begin{proof}[Proof of \eqref{t_birational}] We proceed by induction on the 
dimension.  Assume that \eqref{t_birational} holds up to dimension $n-1$.  By assumption
there is a positive integer $r_{n-1}$ such that $\phi_{rK_W}$ is birational, for any
smooth variety $W$ of general type of dimension at most $n-1$, and for any $r\geq
r_{n-1}$.

Now if $x$ is a very general point of $X$, and $V$ is any subvariety containing $x$, then
any resolution $W$ of $V$ is of general type.  Hence if $\vol(K_X)\geq 1$, then we are
done by \eqref{p_effective}.  Thus we may assume that $\vol(K_X)<1$.  By \eqref{c_lower},
it follows that $X$ is birational to a subvariety of projective space of bounded degree.
But then the set of all varieties of general type such that $\vol(K_X)<1$ is birationally
bounded, and we are done by \eqref{l_bb}.
\end{proof}

\begin{proof}[Proof of \eqref{c_1}] By \eqref{t_birational} and \eqref{l_low} it follows
that any smooth variety of general type of dimension $n$, where $\vol(K_X)\leq M$, is
birational to a subvariety of projective space of degree at most $(r_n)^nM$, and 
the result follows.  
\end{proof}

\begin{proof}[Proof of \eqref{c_2}] If $\vol(K_X)\geq 1$ there is nothing to prove,
otherwise we apply \eqref{c_1} and \eqref{l_bb}.  \end{proof}

\begin{proof}[Proof of \eqref{c_3}] This is an immediate consequence of the main result of 
\cite{Maehara83}.
\end{proof}

\bibliographystyle{hamsplain} 

\bibliography{/home/mckernan/Jewel/Tex/math}

\end{document}
\end